\documentclass{amsart}
\usepackage{amsmath}
\usepackage{paralist}
\usepackage[misc]{ifsym}
\usepackage{epsfig} %For pictures: screened artwork should be set up with an 85 or 100 line screen
\usepackage{epstopdf} %This is to transfer .eps figure to .pdf figure; please compile your article using PDFLaTex or PDFTeXify.
\usepackage[colorlinks=true]{hyperref}
\hypersetup{urlcolor=blue, citecolor=red}
\allowdisplaybreaks

\usepackage{multirow}	% for multirow
\usepackage{makecell}
\usepackage{graphicx}
\usepackage{cleveref}
\usepackage{algorithm}
\usepackage[noend]{algpseudocode}
\usepackage{url}
\algnewcommand{\Inputs}[1]{%
	\State \textbf{Inputs:}
	\Statex \hspace*{\algorithmicindent}\parbox[t]{.8\linewidth}{\raggedright #1}
}
\algnewcommand{\Initialize}[1]{%
	\State \textbf{Initialize:}
	\Statex \hspace*{\algorithmicindent}\parbox[t]{.8\linewidth}{\raggedright #1}
}

\newcommand{\C}{\mathbb{C}}
\newcommand{\R}{\mathbb{R}}

\newcommand{\FW}{\texttt{FW}}
\newcommand{\OF}{\texttt{OF}}
\newcommand{\DT}{\texttt{DT}}
\newcommand{\OFV}{\texttt{Cheat-OF}}
\newcommand{\OFop}{M}

\newcommand{\RE}{\mathrm{re}}
\newcommand{\IM}{\mathrm{im}}

\newtheorem{theorem}{Theorem}[section]

\theoremstyle{definition}

\newtheorem{remark}[theorem]{Remark}

\title[Complex optical flow in dynamic MRI]
{Complex extension of optical flow and its practical evaluation for undersampled \\ dynamic MRI} % Only the first word and proper nouns should be capitalized.

\author[Matthias J. Ehrhardt and Marco Mauritz]{}


\thanks{$^*$Corresponding author: Marco Mauritz}

\begin{document}
\maketitle

% Enter the first author's name and email address; email addresses are required for each author.
% Use footnote notations to indicate address and affiliations with commas between numbers if more than one address applies; see below for examples.
\centerline{\scshape
Matthias J. Ehrhardt $^{{\href{mailto:m.ehrhardt@bath.ac.uk}{\textrm{\Letter}}}1}$
and Marco Mauritz$^{{\href{mailto:marco.mauritz@uni-muenster.de}{\textrm{\Letter}}}*2}$}

\medskip

{\footnotesize
% Enter the full affiliation and country name:
% Do not insert commas or periods at the end of lines.
 \centerline{$^1$University of Bath, United Kingdom}
} % Do not forget to end {\footnotesize with the sign }

\medskip

{\footnotesize
 % Enter the full affiliation and country name:
 \centerline{$^2$University of Münster, Germany}
}

\bigskip

% The name of the handling editor will be entered by AIMS production staff.
% "Communicated by Handling Editor" is not needed for special issue.
 \centerline{(Communicated by Jennifer Mueller)}

%%%%%%%%%%%%%%%%%%%%%%%%%%%%%%%%%%%%%%%%%%%%%%%%%%%%%%%
%             5. ABSTRACT
%%%%%%%%%%%%%%%%%%%%%%%%%%%%%%%%%%%%%%%%%%%%%%%%%%%%%%%

\begin{abstract}
Reconstructing high-quality images from undersampled dynamic MRI data is a challenging task and important for the success of this imaging modality. To remedy the naturally occurring artifacts due to measurement undersampling, one can incorporate a motion model into the reconstruction so that information can propagate across time frames. Current models for MRI imaging are using the optical flow equation. However, they are based on real-valued images. Here, we generalise the optical flow equation to complex-valued images and demonstrate, based on two real cardiac MRI datasets, that the new model is capable of improving image quality.
\end{abstract}

%%%%%%%%%%%%%%%%%%%%%%%%%%%%%%%%%%%%%%%%%%%%%%%%%%%%%%
%                   6. BODY
%%%%%%%%%%%%%%%%%%%%%%%%%%%%%%%%%%%%%%%%%%%%%%%%%%%%%%

\section{Introduction}
\begin{figure}
	\centering
	\includegraphics[width=0.75\linewidth]{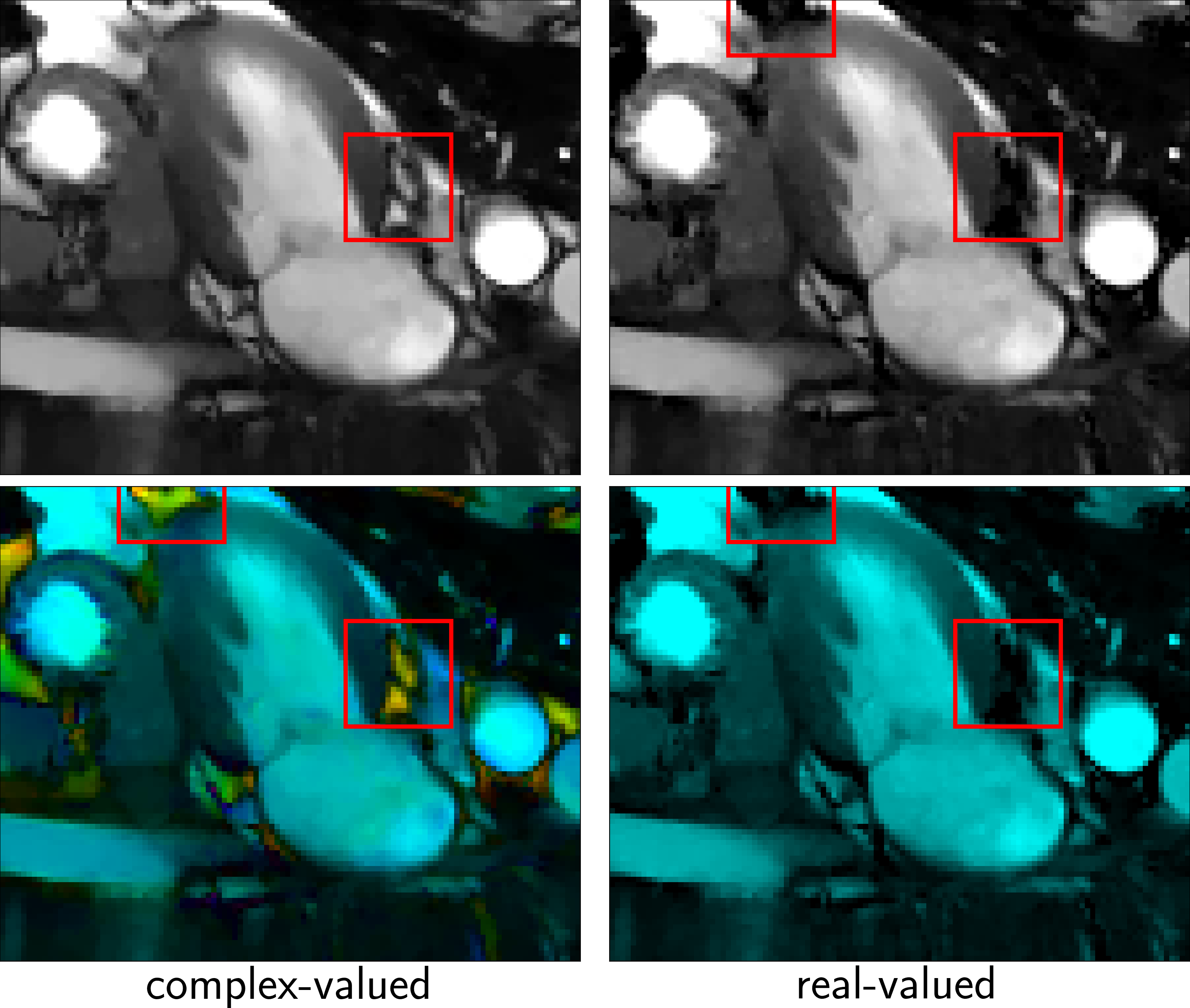}
	\caption{Comparison between real-valued and complex-valued reconstructions in MRI. Left: Complex-valued reconstruction. Right: Real-valued reconstruction. Top row: Magnitude representation. The bottom row shows the same images but as a complex-valued image: The hue represents the normalised phase and the brightness represents the magnitude. The red boxes point to some missing structures in the real-valued image.}
	\label{fig:RealValuedRecon}
\end{figure}

Magnetic resonance imaging (MRI) is a powerful medical imaging modality that allows for noninvasive visualisations of the human body. In particular, for dynamic MRI, there is a variety of applications, e.g.~ cardiac imaging \cite{MRIExample1, MRIExample2}, dynamic contrast-enhanced magnetic resonance imaging \cite{MRIExample3, MRIExample4} and flow imaging \cite{MRIExample5}. However, the imaging times in MRI are inherently slow, leading to motion artifacts that degrade image quality. This issue can be addressed from a software perspective where one aims for acquisition speedups by decreasing the number of measurements and using mathematical modelling to improve image quality. In this article, we follow this approach and solve the ill-posed inverse problem $y_t = A_t\rho^\dagger_t$ (the index $t$ denotes the $t$-th motion state of the considered subject; the MRI forward operator $A_t$ includes subsampling of measurements and hence may depend on $t$) of undersampled dynamic MRI by considering the variational problem
\begin{align*}
	\min_{\rho,v} \sum_t\|A_t\rho_t - y_t\|^2+\mathcal{R}(\rho,v),
\end{align*}
where $y_t$ is a noisy measurement and $\mathcal{R}$ is a regularisation term that incorporates prior knowledge about the images into the reconstruction. More precisely, we aim to jointly reconstruct the temporal image sequence $\rho$ and the velocity field $v$ associated with the temporal evolution. To this end, we model the underlying motion using the optical flow equation \cite{jointDenoising_Burger}. This was already done for MRI and real-valued images \cite{DynMRI_Angelica}. Here, we generalise this idea by explicitly taking into account that images are complex-valued, i.e.~ we formulate an optical flow equation for complex-valued images. This is important since omitting the phase may lead to reconstruction artifacts. This problem can be seen in \cref{fig:RealValuedRecon}. The figure shows reconstructions from fully sampled MRI measurements. The images on the right show a nonnegative real-valued reconstruction where some structures, marked by red boxes, of the complex-valued reconstructions on the left are missing. The top row shows the magnitude of the reconstructions, while the bottom row depicts complex values: The hue represents the normalised phase and the brightness represents the magnitude. In particular, one can see that structures where the phase is different are missing in the real-valued images.

\medskip
\subsection{Related work}
Dynamic inverse problems in general have been studied extensively. An overview of the topic can be found in \cite{DynInvProbOverview_Carola, Schuster_2018}.
The particular topic of undersampled (dynamic) MRI has also been widely studied and many different approaches have been proposed to enhance reconstruction quality. Many works leveraged compressed sensing making use of the fact that MRI images can be sparsely represented in suitable domains \cite{LustigSparseMRI2007, CSSense2009, DynamicCSSense2013}. In dynamic imaging, the temporal relation between motion states should be taken into account to improve the results. One approach is to decompose the (temporal) sequence of images into a low-rank part (background between time frames) and a sparse part
(modelling dynamics between frames) \cite{gao2012compressed, OtazoSparseMatrixDecomposition2015} . This way, the image sequence is represented suitably for the undersampled MRI problem.
Another approach is taken in \cite{MRIRecon_Christoph}, where the motion between different motion states is estimated based on the frame-wise undersampled reconstruction of each state. Subsequently, this motion information is used to get a final reconstruction incorporating all motion states into a single reconstruction.
Joint reconstruction and motion estimation approaches are also promising. This was investigated in \cite{DynMRI_Angelica, DynMRI_linearOF} where the
underlying motion was modelled using (locally affine) optical flow.
Besides classical image reconstruction, learning-based approaches have been considered in dynamic MRI. Those aim at learning spatiotemporal dependencies within the image sequence, and/or mimic the iterative nature of classical reconstruction algorithms \cite{DL_MRI2, DL_MRI1, DL_MRI3}.

\medskip
Besides MRI, dynamic inverse problems arise in other fields of imaging. For instance, in dynamic PET imaging the continuity equation was used as the motion model to describe temporal changes between time points \cite{DynamicPET_Schmitzer}. In tomography, an approach based on a template and its temporal deformation was taken \cite{Chen_2019}.

\medskip
When it comes to working with complex-valued images, two approaches are conceivable: First, one can parameterise the images using real and imaginary part. This is the standard approach that is basically always followed. However, in MRI imaging often the magnitude of the image is of most interest and one would like to pose prior assumptions on the magnitude. This is difficult to incorporate using real and imaginary part because calculating the magnitude is a nonlinear operation, hence the overall regularisation would involve a nonlinear operator. Therefore, one can work with a magnitude and phase parameterisation of the complex-valued images. This makes regularising the magnitude easier but comes at the cost that now the MRI forward operator is nonlinear. This was considered in \cite{Fessler2004, Valkonen_2014}.

\newpage 
\section{Methods}

\subsection{MRI modelling}
Solving the inverse problem in MRI imaging in its simplest form leads to inverting the discrete Fourier transform. Given enough measurements (i.e.~ enough coverage of Fourier space/k-space) this problem is well-posed by the stable invertibility of the Fourier transform. In practice, however, multiple receiver coils (here are the measurements detected) and undersampling lead to an ill-posed inverse problem. The MRI forward operator maps a complex-valued image $\rho_t$ (which in our case is a time slice of a temporal sequence of images) to some complex space $Y_t$, where $Y_t$ depends on the undersampling factor and the number of receiver coils. For instance, in a continuous setting one could choose real and imaginary parts of the image $\rho$ as elements of $L^p((0,T);BV(\Omega))$, hence  the respective time marginals are in $BV(\Omega)$ \cite{jointDenoising_Burger}, where
\begin{align*}
BV(\Omega)=\left\{u\in\mathrm{L}^1(\Omega) \ \middle| \ \sup_{\phi\in C_c^1(\Omega;\R^2), ||\phi||_{\mathrm{L}^\infty}\le1}\int_\Omega u\mathrm{div}\phi\ \mathrm{d}x<\infty\right\}
\end{align*}
is the space of functions of bounded variation \cite{Ambrosio_FunctionsOfBoundedVariation} on $\Omega\subset\R^2$. Similar spaces with more regularity are also possible. In a discrete setting (which will be considered in this article) we take $\rho\in X=\C^{N_t\times N_x\times N_y}$, where $N_t$, $N_x$ and $N_y$ denote the number of dimensions in the temporal, $x$ and $y$ direction, respectively. Therefore, it is $\rho_t\in X_t=\C^{N_x\times N_y}$ for each time slice. We next introduce the forward operator describing the measurements. Denoting by $N_c$ the number of receiver coils, we have $Y_t=\C^{N_c\times N_{k,t}}$, $t=1,\ldots,N_t$, as the measurement space. Here, $N_{k,t}$ denotes the number of measured k-space points at time $t$ after the undersampling has been applied. With the Fourier operator $F$, the coil sensitivity operators $C_{t,i}$ and the undersampling operators $S_{t}$ (see \cref{sec:undersampling} for details) we can find for the $i$-th component (i.e.~ coil) of the MRI forward operator $A_t$ at time point $t$
\begin{align}\label{eq:MRIOperator}
	A_{t}\colon X_t\to Y_t, \quad (A_t\rho_t)_i = S_{t}FC_{t,i}\rho_t.
\end{align}
Hence, we model MRI measurements as
\begin{align*}
	y_t = A_{t}\rho^\dagger_t+\varepsilon_t
\end{align*}
where $\rho^\dagger_t\in X_t$ is the ground truth image that we want to reconstruct and $\varepsilon_t$ is some additive noise.

\subsection{Modelling motion for complex-valued images}\label{sec:ComplexOF}
In this work, we use the optical flow equation as a motion model to regularise the dynamic inverse problem. The optical flow equation is often used in computer vision and has proven valuable for 2D images \cite{DynMRI_Angelica, Burger_2017, jointDenoising_Burger}. Its defining property is constant image intensity along a trajectory with $\frac{\mathrm{d}}{\mathrm{d}t}\rho=v$ \cite{DynInvProbOverview_Carola}. Thus, the underlying PDE reads
\begin{align*}
	\partial_t\rho + \langle v,\nabla\rho\rangle_\R=0
\end{align*}
and during the reconstruction, both the image $\rho$ as well as the velocity field $v$ have to be reconstructed. Note, that since velocity and density are multiplied by each other, using the optical flow equation makes the problem nonconvex. To be able to apply the optical flow equation in MRI imaging (where complex-valued images are reconstructed), we propose to generalise this equation by allowing both the image $\rho$ and the velocity $v$ to be complex-valued. The inner product is therefore replaced by the standard inner product for complex vector spaces. This leads us to the model
\begin{align}\label{eq:ComplexOFEquation}
	0=\partial_t(\rho_\RE+i\rho_\IM)+\langle v_\RE+iv_\IM,\nabla(\rho_\RE+i\rho_\IM) \rangle_\C\\
	\Longleftrightarrow\begin{cases}
		0=\partial_t\rho_\RE+\langle v_\RE,\nabla\rho_\RE\rangle_\R+\langle v_\IM,\nabla\rho_\IM\rangle_\R \nonumber\\
		0=\partial_t\rho_\IM+\langle v_\IM,\nabla\rho_\RE\rangle_\R-\langle v_\RE,\nabla\rho_\IM\rangle_\R\nonumber.
	\end{cases}
\end{align}
for $\rho = \rho_\RE +i\rho_\IM\in X$ and $v = v_\RE +i v_\IM\in V^2$ (note, that the velocity is a 2D vector). Here, $V$ is a suitable space for the velocities such that the derivatives needed for the priors are well-defined. In the discrete setting, we may choose $V=X$.
Associated with the above two equations we introduce the optical flow operator
\begin{align}\label{eq:OFOperator}
	\OFop(\rho, v) = \begin{pmatrix}
		\partial_t\rho_\RE+\langle v_\RE,\nabla\rho_\RE\rangle_\R+\langle v_\IM,\nabla\rho_\IM\rangle_\R\\
		\partial_t\rho_\IM+\langle v_\IM,\nabla\rho_\RE\rangle_\R- \langle v_\RE,\nabla\rho_\IM\rangle_\R
	\end{pmatrix}.
\end{align}
For our discrete space $X=\C^{N_t\times N_x\times N_y}$ the temporal derivative $\partial_t\colon X\to X$ and the inner product $\langle\cdot,\cdot\rangle_\R\colon X\times X\to X$ are to be understood pointwise.

\medskip
We use the same discrete versions of the differential operators as in \cite{jointDenoising_Burger}. This means we use forward differences in time and central differences in space with Neumann boundary conditions for the differential operators involved in the optical flow equation. For the differential operators of other parts of the optimisation function we use forward differences with Neumann boundary conditions (see \cite{jointDenoising_Burger} for more details).

\begin{remark}[Choices on optical flow model]
	There are other possible generalisations of the optical flow equation to complex-valued images. For instance, one could introduce separate velocities for real and imaginary part, hence decoupling both components. Also, a single velocity for both the real and imaginary part would be an option. Finally, only one of the two parts could be regularised by the optical flow equation and the other one evolves freely. However, these variants do not explicitly model interchange between real and imaginary part. To be able to model such an interchange, we used the above approach.
\end{remark}

\begin{remark}[Alternative choice for motion model]
	An alternative motion model would be to use the continuity equation in density and momentum variables instead of density and velocity variables (see, for example, \cite{BeBr00}). This way, the motion model does not incorporate a nonconvexity into the optimisation problem. Unfortunately, finding meaningful regularisers for the momentum variables is more challenging than for the velocities. For instance, an object moving along the $x$ axis with constant velocity will have a zero gradient of the corresponding velocities. The momentum, however, is nonconstant as it also changes when the mass of the underlying object changes. Therefore, we decided to use the optical flow equation.
\end{remark}

\section{Optimisation and implementation}
Based on \cite{DynMRI_Angelica, jointDenoising_Burger} we minimise
\begin{align}\label{eq:optimisationFunction}
	F(\rho, v) = \sum_t||A_t\rho_t-y_t||^2 +
	\alpha_{1}\mathcal{R}_1(\rho)+
	\alpha_{2}\mathcal{R}_2(v)+
	\alpha_{3}&\mathcal{R}_3(\OFop(\rho,v))
\end{align}
with respect to $\rho=\rho_\RE+i\rho_\IM,\ v = (v_{x,\RE},v_{y,\RE})+i (v_{x,\IM},v_{y,\IM})$.
Here, $y_t$ is the measured MRI signal at time point $t$, and $A_t$ and $\OFop$ are operators related to the MRI measurements and the optical flow equation, respectively (see \eqref{eq:MRIOperator} and \eqref{eq:OFOperator}). The regularisers for the variables $\rho$ and $v$ are chosen to promote smoothness (the ``level of smoothness'' can be controlled via the $\varepsilon_i$ parameters)
\begin{align*}
	&\mathcal{R}_1(\rho) = \mathcal{H}_{\varepsilon_1}(\nabla\rho_\RE) + \mathcal{H}_{\varepsilon_1}(\nabla\rho_\IM)\\
	&\mathcal{R}_2(v) = \mathcal{H}_{\varepsilon_2}(\nabla v_{x, \RE}) + \mathcal{H}_{\varepsilon_2}(\nabla v_{y,\RE})+ \mathcal{H}_{\varepsilon_2}(\nabla v_{x,\IM})+ \mathcal{H}_{\varepsilon_2}(\nabla v_{y,\IM}),
\end{align*}
where $\nabla$ denotes the spatial gradient and $\mathcal{H}_\varepsilon$ is the Huber loss with parameter $\varepsilon$. Defining for $x\in\R^2$
\begin{align*}
	h_{\varepsilon}(x) =
	\begin{cases}
		\frac{\|x\|^2_2}{2\varepsilon} \ \ &\text{if }\|x\|_2\le\varepsilon,\\
		\|x\|_2 - \frac{\varepsilon}{2} \ \ &\text{else},
	\end{cases}
\end{align*}
the Huber loss is defined as $\mathcal{H}_\varepsilon(x)=\sum_{i=1}^{N}h_\varepsilon(x_i)$ for $x\in (\R^2)^N$, and it can be seen as a smooth approximation of the norm $||\cdot||_2$ \cite{HuberLoss}. Finally, we use a soft constraint
\begin{align*}
	\mathcal{R}_3(\OFop(\rho,v)) = \mathcal{H}_{\varepsilon_3}(\OFop(\rho,v))
\end{align*}
since we cannot expect the optical flow equation to be satisfied perfectly. Compared to the classical $ \mathrm{L}^1$ loss the Huber loss is differentiable, and the parameter $\varepsilon$ controls the extent of smoothness within the Huber loss, i.e.\ it controls whether the Huber loss is more similar to the norm $\|\cdot\|_2$ or the squared norm $\|\cdot\|_2^2$. Tuning this parameter allows us to adapt to the considered problem and find a better model.\medskip

% The weight parameters
% $\alpha_{\mathrm{OF}}$, $\alpha_\rho$ and $\alpha_v$ as well as the Huber parameters $\varepsilon_{\mathrm{OF}}$, $\varepsilon_\rho$ and $\varepsilon_v$ are tunable hyperparameter that will be optimised for each considered dataset using bayesian optimisation.

The considered optimisation problem \eqref{eq:optimisationFunction} is nonconvex due to the optical flow equation. To deal with the nonconvexity we follow the same strategy as in \cite{DynMRI_Angelica, jointDenoising_Burger} and apply a block coordinate descent (see \cref{algo:Optimisation}). This means we iteratively optimise the minimisation function only for $\rho$ or $v$ while keeping the other variable fixed. Convergence of this method is given by \cite[Thm. 4.1]{BlockCoordinateDescent} if we additionally assume the norms of $\rho$ and $v$ to be bounded. The bound can be arbitrarily large and is only needed for theoretical reasons: It was never reached in the numerical examples. This is a natural assumption since mass and velocity have to stay within certain bounds in reality. Unfortunately, no convergence rates are established in \cite{BlockCoordinateDescent}, and one has to decide to what degree of accuracy the subproblems are being solved. We empirically observed that the overall algorithm produces a converging sequence.\medskip

Each of the subproblems is solved using the FISTA algorithm \cite{FISTAAlgorithm} (applied to $f+g$ with $f$ smooth and $g=0$, see \cref{algo:FISTA}) with early stopping (iterations where stopped if the relative $\mathrm{L}^2$ error was below $\delta=10^{-5}$). Furthermore, motivated by the coarse-to-fine approach in image registration, we smoothed the density $\rho$ and the velocity $v$ after each sub-optimisation using a Gaussian smoothing with standard deviation $\sigma / {i_\mathrm{outer}}$. Here, $\sigma$ is again a tunable parameter, and $i_\mathrm{outer}\in\{1,\ldots,n_\mathrm{outer}\}$ denotes the current outer iteration, i.e.\ the smoothing decays with growing iterations. For the maximal number of iterations we used $n_\mathrm{outer}=200$, $n_{\rho}=1400$ and $n_{v}=3200$.
The optimisation algorithm was implemented in Python using SIRF \cite{SIRF_paper} for reading the MR data and applying the associated measurement operators, and CIL \cite{CIL_paper} to implement the actual optimisation algorithm. The code is available on \url{https://github.com/mjmaur/dynamic_MRI}.

\begin{algorithm}
	\caption{FISTA with early stopping for minimizing smooth function $f$}
	\label{algo:FISTA}
	\begin{algorithmic}[1]
		\Inputs{smooth function $f$ with Lipschitz constant $L$, initialisation $x_\mathrm{init}$\\
			Algorithm parameters $n$, $\delta$}
		\Initialize{$x^0\gets x_\mathrm{init}$, $\hat x^1\gets x_\mathrm{init}$, $t^1\gets 1$}
		\For{$j=1$ to $n$}
		\State $x^{j}\gets$ $\hat x^{j}-\frac{1}{L}\nabla f(\hat x^{j})$
		\State $t^{j+1}\gets\frac{1+\sqrt{1+4(t^j)^2}}{2}$
		\State $\hat x^{j+1} \gets x^{j}+(\frac{t^j-1}{t^{j+1}})(x^{j} - x^{j-1})$
		\If{$\|x^{j}-x^{j-1}\|_2 / \| x^{j-1}\|_2 < \delta$}
		\State break
		\EndIf
		\EndFor
		\State return $x^j$
	\end{algorithmic}
\end{algorithm}

\begin{algorithm}
	\caption{Optimisation of $F(\rho,v)$}
	\label{algo:Optimisation}
	\begin{algorithmic}[1]
		\Inputs{Measurement $y$\\
			Model parameters $\alpha_1,\alpha_2,\alpha_3,\varepsilon_1,\varepsilon_2,\varepsilon_3$\\
			Algorithm parameters $\sigma,n_\mathrm{outer}, n_{\rho}, n_{v}$, $\delta$}
		\Initialize{$v^1 \gets 0$, $\hat v^1\gets 0$, $\rho^1\gets 0$, $\hat\rho^1\gets 0$}
		
		\For{$i = 1$ to $n_\mathrm{outer}$}
		\State $\rho^{i+1}\gets \text{FISTA}(\rho\mapsto F(\rho,\hat v^i),\hat\rho^i)$ with $n_\rho$ iterations
		\State $\hat\rho^{i+1}\gets\mathrm{smooth}_{\sigma/i}(\rho^{i+1})$
		\State $v^{i+1}\gets \text{FISTA}(v\mapsto F(\hat\rho^{i+1},v),\hat v^i)$ with $n_v$ iterations
		\State $\hat v^{i+1}\gets\mathrm{smooth}_{\sigma/i}(v^{i+1})$
		\If{$\frac12 \|v^{i+1}-v^{i}\|_2 / \| v^{i}\|_2 + \frac12 \|\rho^{i+1}-\rho^{i}\|_2 / \| \rho^{i}\|_2< \delta$}
		\State break
		\EndIf
		\EndFor
		\State return $\rho^{i+1}, v^{i+1}$
	\end{algorithmic}
\end{algorithm}

\subsection{Models for comparison}
We compare our proposed method \eqref{eq:optimisationFunction} (we call this model \OF~ and label corresponding reconstructions \OF, too) to three slightly different variants:
\begin{itemize}
	\item[(a)] Frame-wise model (\FW): For each time point a separate reconstruction is computed and no correlation between time points is assumed. This corresponds to optimising \eqref{eq:optimisationFunction} for the density $\rho$ only and setting $\alpha_2$ and $\alpha_3$ to zero. Thus, we minimise
	\begin{align*}\label{eq:optimisationFunctionFW}
		F(\rho) = \sum_t||A_t\rho_t-y_t||^2 +
		\alpha_{1}\mathcal{R}_1(\rho)\nonumber
	\end{align*}
	with respect to $\rho=\rho_\RE+i\rho_\IM$ applying \cref{algo:FISTA} with $n_\rho$ iterations, where
	\begin{align*}
		\mathcal{R}_1(\rho) = \mathcal{H}_{\varepsilon_1}(\nabla\rho_\RE) + \mathcal{H}_{\varepsilon_1}(\nabla\rho_\IM)
	\end{align*}
	as above.
	\item[(b)] Time derivative model (\DT): The term containing the optical flow equation is replaced by a term containing time derivatives only. This means we fix the velocities to zero and optimise \eqref{eq:optimisationFunction} for the density $\rho$ only. Thus, we minimise
	\begin{align*}%\label{eq:optimisationFunctionDT}
		F(\rho) = \sum_t||A_t\rho_t-y_t||^2 + \alpha_{1}\mathcal{R}_1(\rho)+
		\alpha_{3}\mathcal{R}_3(\OFop(\rho, 0))\nonumber
	\end{align*}
	with respect to $\rho_\RE+i\rho_\IM$ applying \cref{algo:FISTA} with $n_\rho$ iterations, where
	\begin{align*}
		&\mathcal{R}_1(\rho) = \mathcal{H}_{\varepsilon_1}(\nabla\rho_\RE) + \mathcal{H}_{\varepsilon_1}(\nabla\rho_\IM)
	\end{align*}
	and
	\begin{align*}
		&\mathcal{R}_3(\OFop(\rho,0)) = \mathcal{H}_{\varepsilon_3}([\partial_t\rho_\RE, \partial_t\rho_\IM])
	\end{align*}
	as above promote spatial and temporal regularity, respectively.\medskip
	
	\item[(c)] Optical flow using ground truth velocities (\OFV): With this model we want to investigate how our method would perform if we had access to the velocities underlying the motion, i.e.\ this model gives us the best possible reconstruction within our framework. To this end, we use our model \eqref{eq:optimisationFunction} with fully sampled data to determine ground truth velocities $v_\mathrm{GT}$. Then, we optimise \eqref{eq:optimisationFunction} for the density $\rho$ only while fixing the velocity to be $v_\mathrm{GT}$. Thus, we minimise
	\begin{align*}%\label{eq:optimisationFunctionOFV}
		F(\rho) = \sum_t||A_t\rho_t-y_t||^2 +
		\alpha_{1}\mathcal{R}_1(\rho)+\alpha_{3}\mathcal{R}_3(\OFop(\rho, v_{\mathrm{GT}}))
	\end{align*}
	with respect to $\rho=\rho_\RE+i\rho_\IM$ applying \cref{algo:FISTA} with $n_\rho$ iterations, where
	\begin{align*}
		&\mathcal{R}_1(\rho) = \mathcal{H}_{\varepsilon_1}(\nabla\rho_\RE) + \mathcal{H}_{\varepsilon_1}(\nabla\rho_\IM)
	\end{align*}
	and
	\begin{align*}
		&\mathcal{R}_3(\OFop(\rho, v_{\mathrm{GT}}))=\mathcal{H}_{\varepsilon_3}(\OFop(\rho,v_{\mathrm{GT}})
	\end{align*}
	as above.
\end{itemize}

\subsection{Metrics}
We consider two metrics to compare the quality of two images: Peak signal-to-noise ratio (abbreviated as PSNR; this metric is related to the mean squared error, and is used as the metric for the hyperparameter optimisation) and structural similarity index measure \cite{SSIM} (abbreviated as SSIM; this is the secondary alternative metric to compare image quality). To compute them, we use the respective functions provided by \textit{scikit-image} where the ``data\_range" parameter is obtained from the ground truth image. Also, we compute the metrics for each time point and compute the final metric as the average over all time points.

\subsection{Hyperparameter optimisation}
Hyperparameter optimisation was carried out using \textit{scikit-optimize} to find good sets of parameters for each model. We chose PSNR (computed on the dynamic mask, see \cref{sec:DynamicMask}) as the metric to be minimised by the Bayesian optimisation. For the Bayesian optimisation we used the Optimizer class with a Gaussian process as the base estimator and expected improvement as the acquisition function. For the model in \eqref{eq:optimisationFunction} and its optimisation we used the following hyperparameters. The regularisation parameters $\alpha_1, \alpha_2,\alpha_3$, the smoothness parameters $\varepsilon_1,\varepsilon_2,\varepsilon_3$ and the algorithm parameter $\sigma$ which gives the initial amount of smoothing of the density and velocity after each outer optimisation step. The number of tested parameters depends on the model. The more model parameters exist, the more combinations of parameters were considered within the Bayesian optimisation. More details can be found in the code.

\section{Data}\label{sec:Data}
We work with two datasets taken from the OCMR dataset \cite{OCMRDataset} containing cardiac cine series (long-axis view). The data consists of multi-coil k-space data acquired with Cartesian sampling on a Siemens MAGNETOM Avanto scanner (1.5T) (see \cite{OCMRDataset} for more details). We use fully sampled data and apply the subsampling manually. A fully sampled reconstruction of both datasets can be found in \cref{fig:GTAndMask}. Based on subject 1, we additionally consider a simulated dataset that resembles a real-world MRI image but satisfies the optical flow equation perfectly. This way, we can assess the model capacities and our optimisation procedure in an idealised setting.

\subsection{Ground truth images}\label{sec:GroundTruth}
To assess the reconstruction quality, we need a suitable ``ground truth'' image to compare the results to. One possibility would be to use the fully sampled measurements to get reconstructions for each time frame separately using the \FW~ model. In doing so, however, one gets time trajectories that are slightly irregular and do not resemble realistic motion (see \cref{fig:GTAndMask} second column on the right) due to, e.g., in-frame motion. Therefore, we use a ``dynamic ground truth'' that results from minimising \eqref{eq:optimisationFunction} with fully sampled data and high dynamic regularisation parameters. This results in reconstructions that are more regular in time (see \cref{fig:GTAndMask} first column on the right). This approach also has the advantage that we have ``ground truth'' velocities available that describe the dynamics within the ground truth image sequence. These can be used to assess the quality of our reconstructed velocities.\medskip

In the discussion section, we will briefly comment on the results being similar when using the frame-wise ``ground-truth''.

\medskip
\subsection{Mask focusing on dynamic parts in images}\label{sec:DynamicMask}
\begin{figure}
	\centering
	\includegraphics[width=1\linewidth]{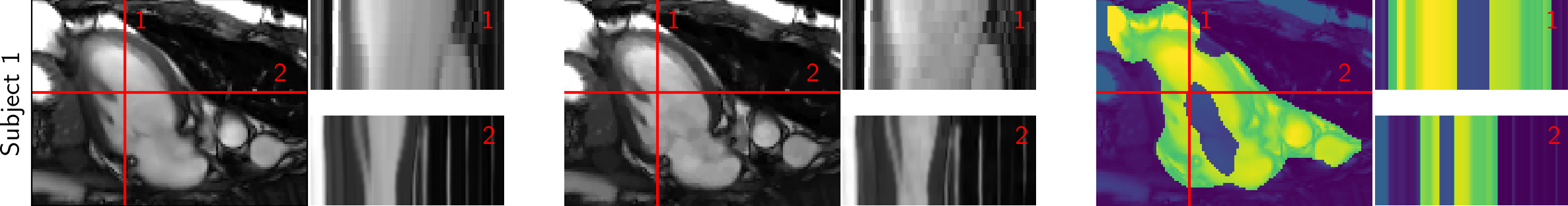}\\
	\vspace{0.5cm}
	\includegraphics[width=1\linewidth]{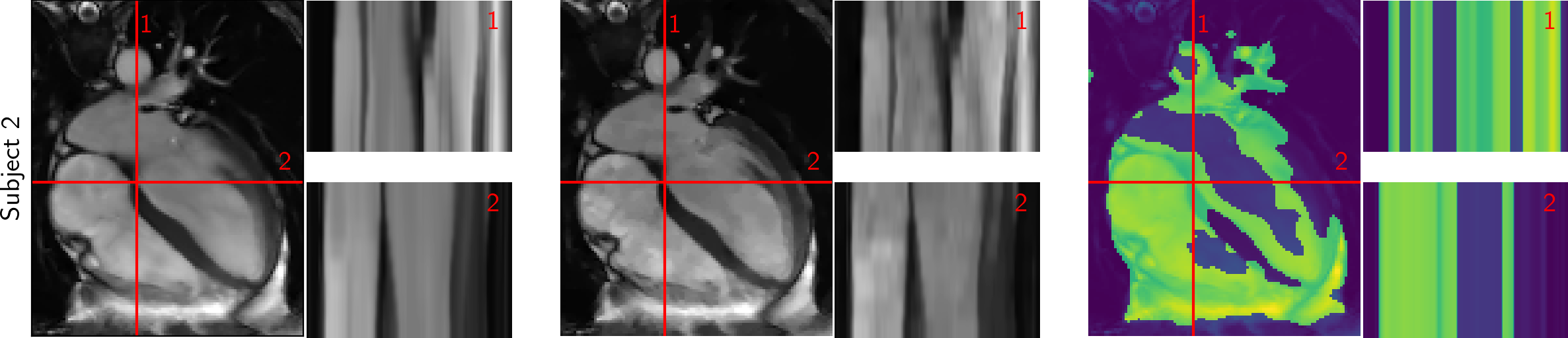}
	\caption{Each row shows a reconstruction for a different data set. The first column shows a dynamic ground truth (obtained as a reconstruction using \eqref{eq:optimisationFunction} and fully-sampled data). The second column shows the same for a frame-wise ground truth (obtained from reconstructing each time frame separately (\texttt{FW} model) using fully-sampled data). The last column shows an overlay of the reconstructions with the dynamic masks. Within each column, the left side shows one time slice of the reconstructions, the right side shows time-space images along the red lines. }
	\label{fig:GTAndMask}
\end{figure}

We are primarily interested in the dynamic parts of the image and want to improve image quality in these regions particularly. Since large parts of the images are basically constant in time, computing loss metrics between ground truth images and reconstructions over the whole image would consequently result in focusing on the wrong parts. Therefore, we computed masks focusing on the dynamic parts of the images (see \cref{fig:GTAndMask} right column), and computed the image metrics on these masks only.

\subsection{Simulated data}
We test our proposed model on simulated data. To this end, we used a dataset from \cite{OCMRDataset} (subject 1 in \cref{fig:GTAndMask}) and solved \eqref{eq:optimisationFunction} for fully sampled data. Then, we smoothed the obtained velocities in space and time and used the resulting velocities to simulate a sequence of images starting with the image at the first time point from solving \eqref{eq:optimisationFunction}. The thus obtained images and velocities are therefore guaranteed to satisfy the optical flow equation and the velocities are smooth enough to be reconstructable. Based on these simulated ground truth images we generated measurements by applying the MRI forward operator $A_t$. Finally, Gaussian noise (zero mean) was added where the standard deviation was chosen such that unregularised reconstructions yielded qualitatively similar images to unregularised reconstructions based on the original measured data.

\subsection{Undersampling}\label{sec:undersampling}
We next describe the used undersampling. The data we are working with has been acquired with a standard Cartesian sampling scheme, and we consider four times undersampling in the phase encoding direction with full sampling of the rows in the readout direction. More precisely, we use the same $x=15$ percent central rows in each time frame and randomly select $25-x$ percent of the outer rows (same number of rows on each side of the central rows) in such a way that the selected outer rows between consecutive time frames do not overlap.

\section{Experiments and results}
\subsection{Simulated data}\label{sec:experimentsSimulatedData}
\begin{figure}
	\centering
	\includegraphics[width=1\linewidth]{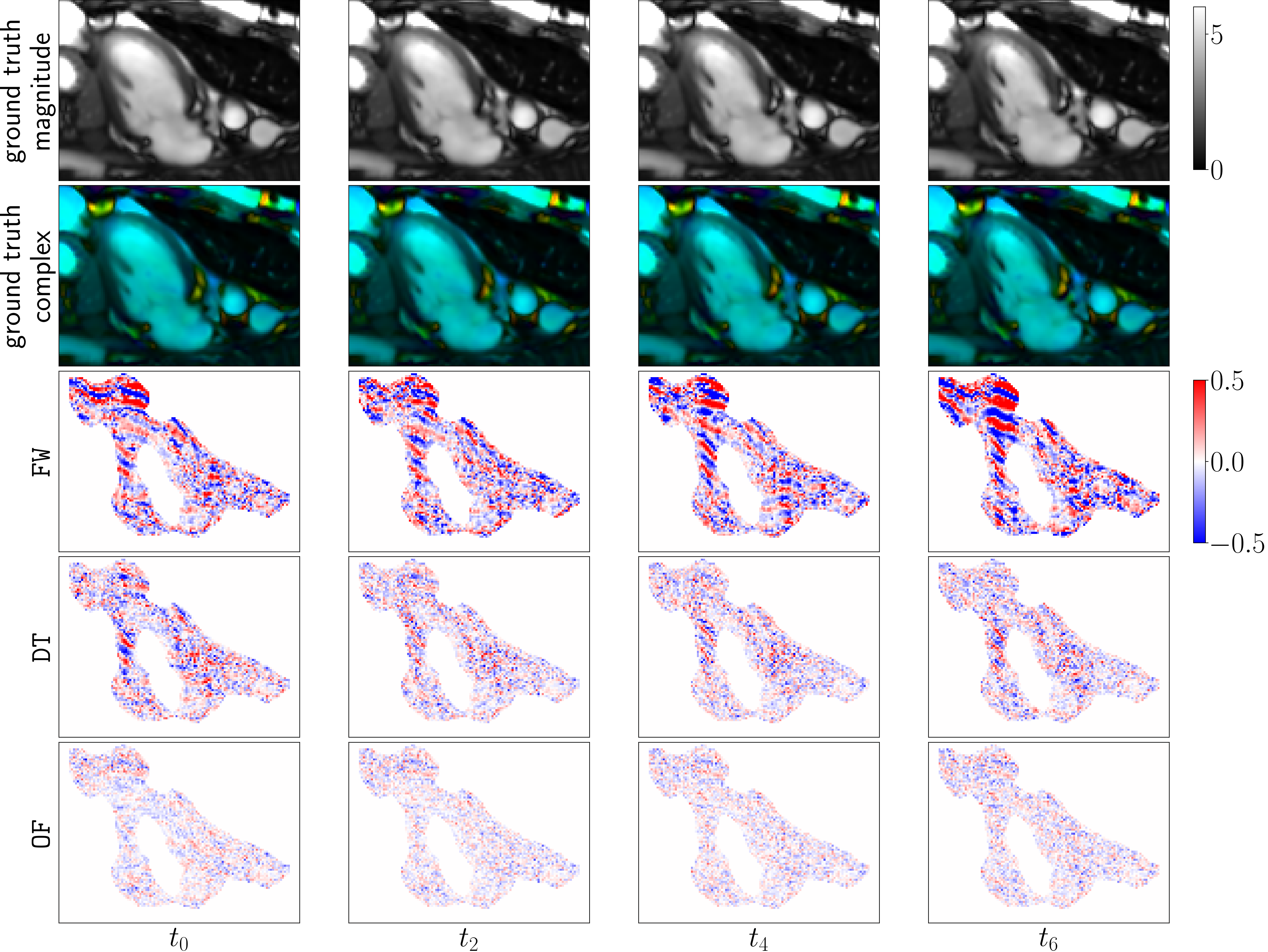}
	\vspace{-0.25cm}
	\caption{Results for the simulated data. The top row shows the evolution (magnitude) of the simulated ground truth images for the time points $t_0$, $t_2$, $t_4$, $t_6$ (out of $8$ time points). The second row shows the same images but as a complex-valued image: The hue represents the normalised phase and the brightness represents the magnitude. The other rows depict the difference between the magnitudes of the reconstructions and the magnitude of the simulated images on the dynamic mask. Rows 3-5: \FW, \DT\, and \OF.}
	\label{fig:SimulationReconstruction}
\end{figure}

The results for the simulated data are summarised in \cref{ResultsTable} and can visually be inspected in \cref{fig:SimulationReconstruction}. The top row in that figure shows the magnitude of the simulated images. In the second row, the same images are represented as complex-valued images where the hue represents the normalised phase and the brightness represents the magnitude. From this row of images, we see that the phase varies spatially, but its temporal changes are small.

The other rows depict the difference between the reconstructions using the respective model and the simulated images (magnitude images are compared). We can clearly see that \DT~ improves over the frame-wise approach \FW, and that \OF~ improves over \DT. The improved reconstruction qualities can also be seen in PSNR values in \cref{ResultsTable} (first column), where \OF~ improves over \FW~ and \DT~ by 10dB and 4dB, respectively. Moreover, we see, that \OFV~, for which ground truth velocities are available in the reconstruction, is basically not improving over \OF. Similar observations can be made by comparing the SSIM values in \cref{ResultsTable}.

\begin{figure}
	\centering
	\includegraphics[width=1\linewidth]{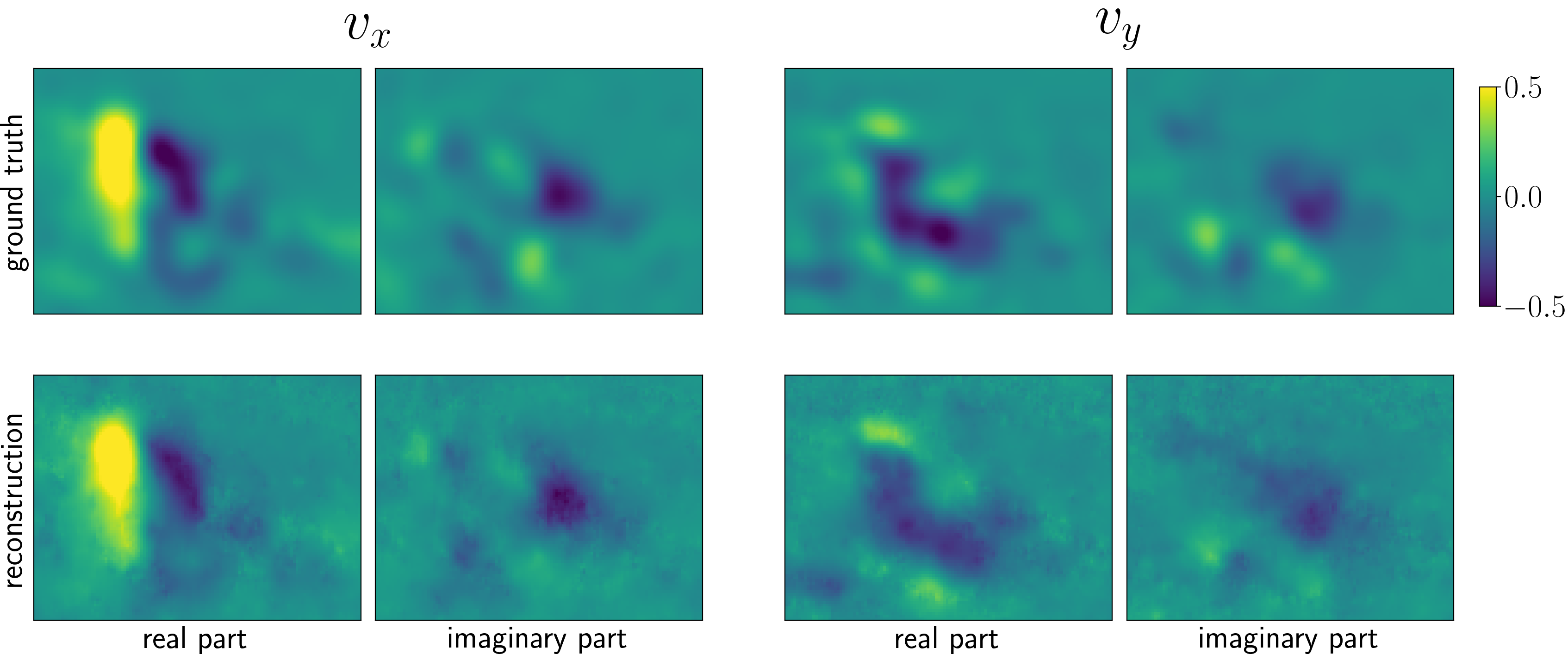}
	\caption{Reconstructed velocities of one time frame for the simulated data using our proposed complex-valued optical flow method. First row: Ground truth velocities in $x$ and $y$ direction. Second row: Reconstructed velocities in $x$ and $y$ direction}
	\label{fig:SimulationVelocities}
\end{figure}

This good reconstruction result of the \OF~ model is backed by the reconstructed velocities depicted in \cref{fig:SimulationVelocities} (bottom row). We see that the general form of the velocities has been recovered compared to the ground truth velocities (top row), hence a good reconstruction quality is expected.

\subsection{Subject 1}
\begin{figure}
	\centering
	\includegraphics[width=0.98\linewidth]{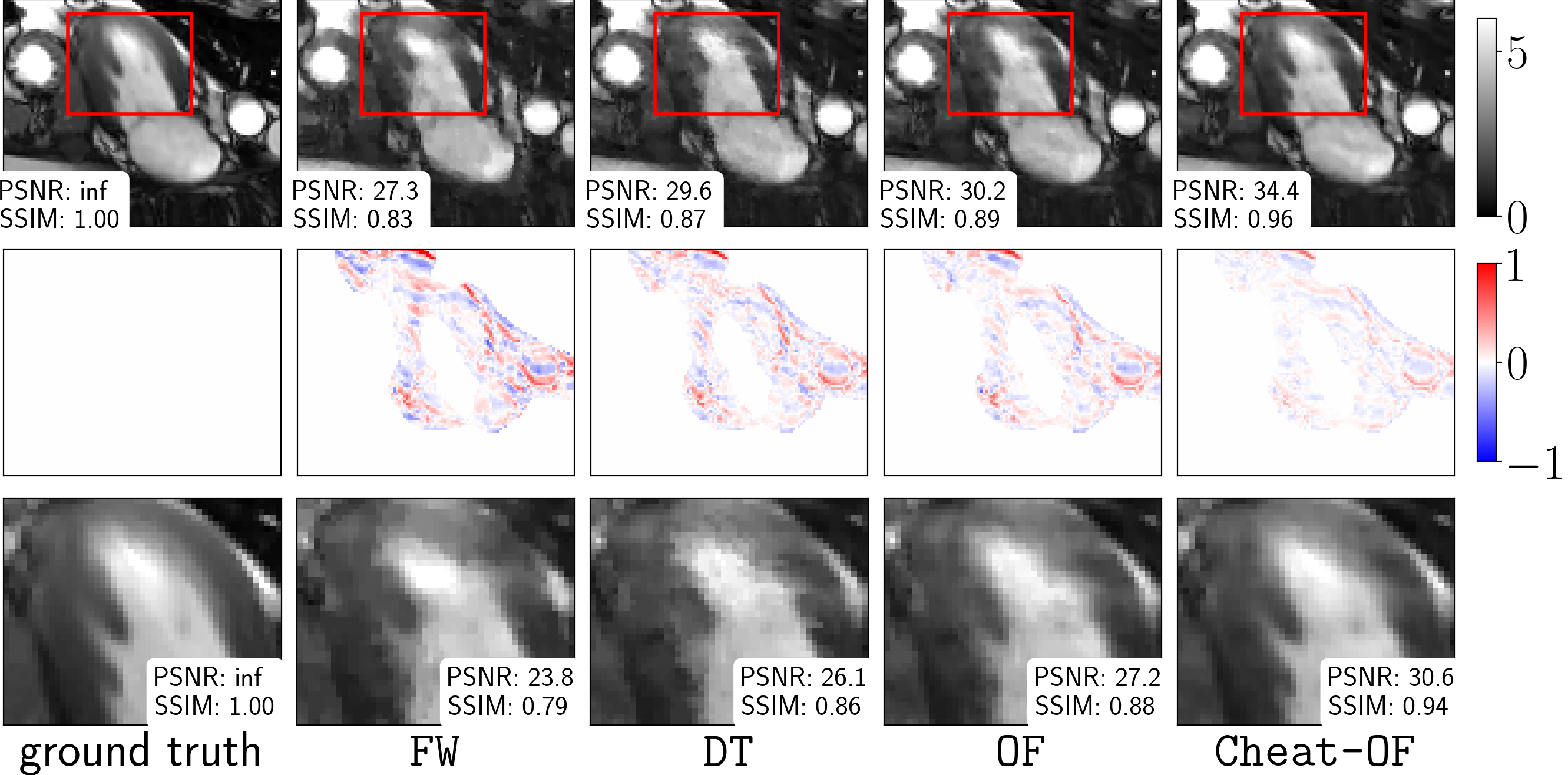}
	\caption{Reconstructions for subject 1. Left to right columns: Ground truth image, \FW, \DT, \OF\ and \OFV. First row shows magnitude of reconstructions. Second row: Difference of reconstruction to ground truth on the mask. Third row: Magnified area and corresponding metrics on magnification.}
	\label{fig:comparison_0001_0002_1}
\end{figure}

\begin{figure}
	\centering
	\includegraphics[width=0.95\linewidth]{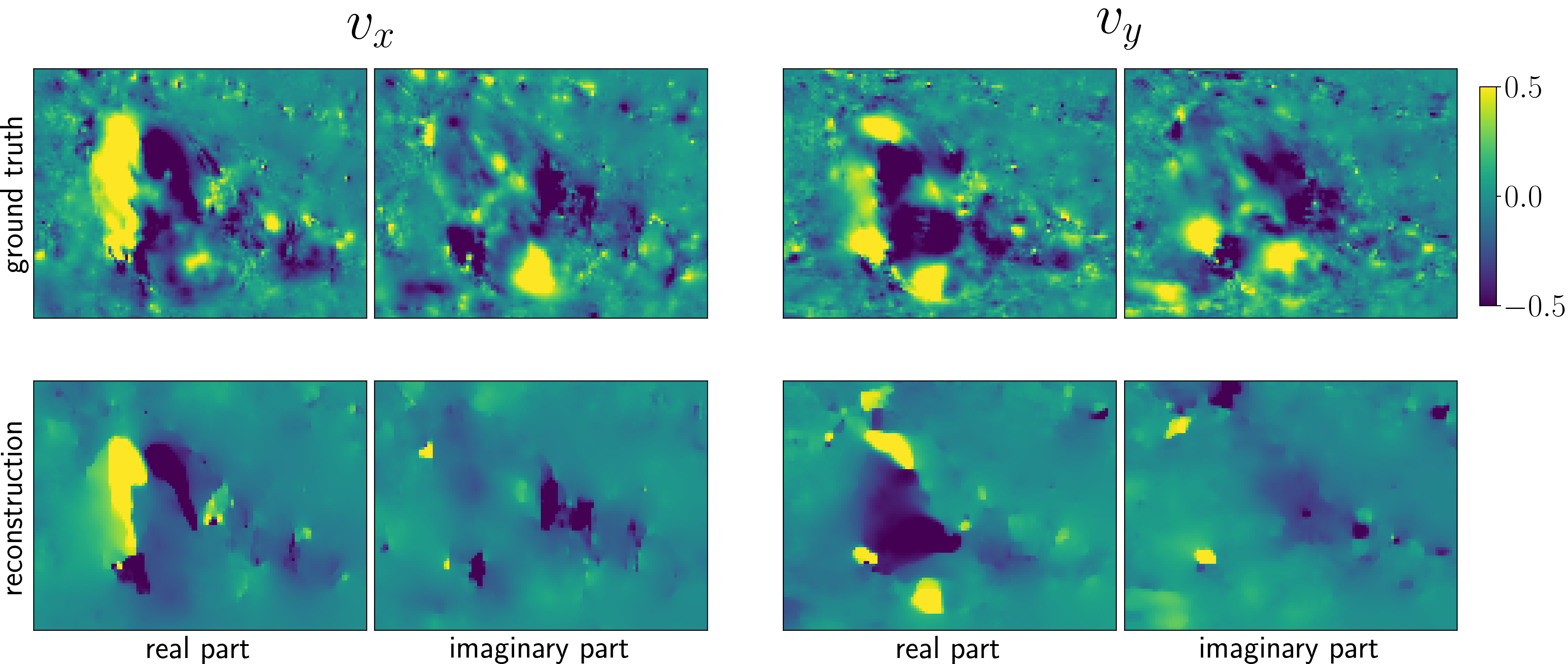}
	\caption{Comparison between ground truth velocities (top row) and reconstructed velocities (bottom row) of subject 1 for a single time frame.}
	\label{fig:comparison_0002_velocity}
\end{figure}

The summarised results of subject 1 can be found in \cref{ResultsTable} and \cref{fig:comparison_0001_0002_1}. We see that our proposed method \OF~ performs better than \FW~ or \DT. The comparison between \OF~ and the frame-wise approach yields PSNR improvements (a similar pattern can be observed in the SSIM values) of almost 3dB, showing that the complex optical flow equation helps the reconstruction by describing the underlying motion. The gains of \OF~ over \DT~ are not very pronounced. However, they improve to over 1dB when looking at finer structures in more detail as can be seen in \cref{fig:comparison_0001_0002_1} in the third row. Further, the results in \cref{ResultsTable} and \cref{fig:comparison_0001_0002_1} also show that the method \OFV~ including ground truth velocities outperforms all the other methods by a lot (which is no surprise since this method basically has access to the fully sampled measurements via the ground truth velocities obtained from a fully sampled reconstruction using \eqref{eq:optimisationFunction}). This improvement can be explained by comparing the reconstructed to the ``ground truth'' velocities in \cref{fig:comparison_0002_velocity}. We see that the reconstruction quality is decent for the real part, whereas the reconstructed imaginary part differs significantly from the ``ground truth'' velocities' imaginary part. Overall, the velocities are not as well reconstructed as for the simulated data, which explains the improved results using \OFV.

\subsection{Subject 2}
\begin{figure}
	\centering
	\includegraphics[width=0.99\linewidth]{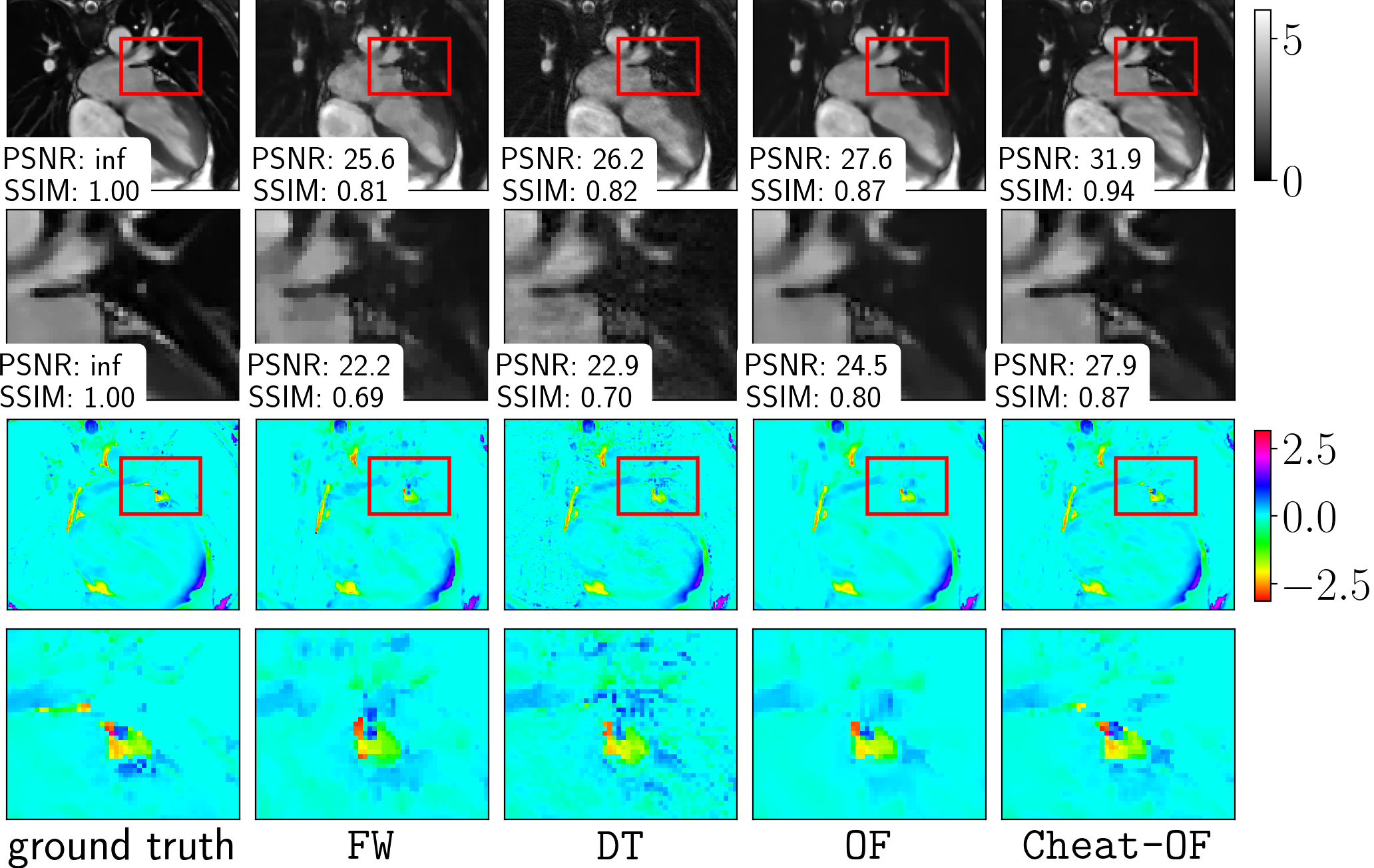}
	\caption{Reconstructions for subject 2. Left to right columns: Ground truth image, \FW, \DT, \OF\ and \OFV. First and second row show magnitude of reconstructions with respective metrics. Third and fourth row show the phase of the reconstructions. Second and last row are the indicated magnified areas.}
	\label{fig:comparison_0001_0002_2}
\end{figure}

The summarised results of subject 2 can be found in \cref{ResultsTable} and \cref{fig:comparison_0001_0002_2}. We see that \OF~ improves PSNR over \FW~ and \DT~ by 2dB and 1.4dB, respectively. Therefore, the proposed optical flow model helps with reconstructing dynamic MRI images. Again, a similar conclusion can be drawn from the SSIM values. Looking at the reconstructed phase (\cref{fig:comparison_0001_0002_2} rows 3 and 4), we see that all methods reconstruct the phase well on an overall level (row 3). However, the magnified area in row 4 shows differences in the reconstructed phases in the detailed structures. Here, \OF~ gives a better result over \FW~ and \DT. This improved phase reconstruction is associated with a better reconstruction of details, as can be seen in the second row. Again, \OFV, as expected, improves over all other methods.

\setlength{\tabcolsep}{1.3pt}

\setlength{\tabcolsep}{4pt}
\renewcommand{\arraystretch}{1.5}
\begin{table}
	\centering
	\scalebox{0.85}{\begin{tabular}{cc||c|c|c||c|c}
		&& \rotatebox[origin=c]{90}{\makecell{\ simulation}} & \rotatebox[origin=c]{90}{\makecell{\ subject 1\\ GT(dyn)}} &  \rotatebox[origin=c]{90}{\makecell{\ subject 2 \\ GT(dyn)}} &
		\rotatebox[origin=c]{90}{\makecell{\ subject 1 \\ GT(FW)}} &
		\rotatebox[origin=c]{90}{\makecell{\ subject 2 \\ GT(FW)}}\\
		\hline
		\hline
		& \FW & $34.1\pm1.2$ & $27.3\pm0.7$ & $25.6\pm0.6$ & $27.7\pm 0.8$  & $26.4\pm 0.9$\\
		& \DT & $40.1\pm0.9$& $29.6\pm0.3$ & $26.2\pm0.3$ & $29.8\pm 0.4$ & $27.5\pm 0.5$\\
		& \OF & $44.1\pm0.6$& $30.2\pm0.5$ & $27.6\pm0.4$ & $30.2\pm 0.5$ & $28.0\pm 0.5$\\
		\cline{2-7}
		\multirow{-4}{*}{\rotatebox[origin=c]{90}{PSNR}}& \OFV & $44.2\pm0.3$& $34.4\pm0.7$ & $32.0\pm0.7$ & $33.1\pm 0.4$ & $ 28.7\pm 1.0$
		\\
		\hline
		\hline
		& \FW &$0.943\pm0.010$& $0.828\pm0.014$ & $0.812\pm0.015$ & $0.847\pm 0.014$ & $0.852\pm 0.018$\\
		& \DT & $0.978\pm0.004$& $0.869\pm0.009$ & $0.819\pm0.012$ & $0.880\pm 0.009$ & $0.869\pm0.011$\\
		& \OF &$0.991\pm0.001$& $0.892\pm0.010$  & $0.873\pm0.012$ & $0.888\pm 0.011$ & $0.880\pm 0.011$\\
		\cline{2-7}
		\multirow{-4}{*}{\rotatebox[origin=c]{90}{SSIM}} & \OFV & $0.992\pm0.001$ & $0.959\pm0.009$ & $0.942\pm0.009$ & $0.934\pm 0.004$ & $0.872\pm 0.018$\\ \\
	\end{tabular}}
	\vspace{-0.4cm}
	\caption{Summary of the results. The columns show the metrics for the respective data sets (with dynamic ground truth GT(dyn)
		and frame-wise ground truth GT(FW)). The displayed values are $(\text{average})\pm(\text{standard deviation})$ taken over the considered time points. The metric for each time point is computed on a dynamic mask (see \cref{sec:DynamicMask}, \cref{fig:GTAndMask}).}
	\label{ResultsTable}
\end{table}


\section{Discussion}
Overall, our experiments demonstrated that explicitly considering the complex-valued optical flow equation \eqref{eq:ComplexOFEquation} is useful for modelling motion between complex-valued images. In particular, our simulated experiment showed large reconstruction improvements. This, of course, is to be expected as the underlying ground truth motion perfectly satisfies the optical flow equation and the corresponding velocities follow the posed prior assumptions. Still, this simplified example shows that the optimisation procedure is suitable for the problem, and meaningful results are obtained. Moreover, we have mentioned already in \cref{sec:experimentsSimulatedData} that \OFV~ is basically not improving over \OF~ in the simulation, suggesting that \OF~ is able to produce nearly optimal reconstructions in case the underlying images satisfy the considered prior assumptions.

Comparing the ground truth velocities of the simulation (\cref{fig:SimulationVelocities}) and of the measured subject 1 (\cref{fig:comparison_0002_velocity}) we see that they are of different quality. In particular, the simulated velocity appears to be much smoother in space, hence they are more in accordance with our prior modelling. On the other hand, the ground truth velocities of subject 1 are much rougher in space which is not reflected in the reconstructed velocity. This of course influences reconstruction quality and partly explains the differences in reconstruction quality compared to the simulation. In future work, it will be interesting to investigate different priors that better resemble the underlying velocities' structure. For instance, one could use Haar wavelets, include regularisation in time, or follow a learning-based approach that can derive very specialised and problem-specific priors.

For our ``ground truth'' images we argued in \cref{sec:GroundTruth} that using ``dynamic ground truth'' images is more appropriate in our setting. However, we repeated the procedures for subjects 1 and 2 with ``frame-wise ground truth'' images (for the simulation both versions of ground truth images coincide) with similar results, where PSNR improvements are slightly less pronounced (2.5dB over \FW~ for subject 1 and 1.6dB over \FW~ for subject 2, see \cref{ResultsTable} right two columns). This suggests that the overall results do not depend on the choice of ``ground truth'' images.

Lastly, we discuss the computational time needed to solve the optimization problems. The scope of this work was to investigate if the proposed complex-valued optical flow equation enhances reconstruction quality. Therefore, we decided to work with conservative hyperparameter choices, i.e. high numbers of inner and outer iterations for \cref{algo:Optimisation} so that the optimization procedure has as little impact as possible on the results. Additionally, the employed software libraries CIL and SIRF are not optimized for computational speed but rather for rapid prototyping. We ran each optimization problem on a single CPU-core which took several hours depending on the model parameters. However, we want to emphasize that a more efficient implementation should not be too complicated to achieve. For example, applying the forward operator \eqref{eq:MRIOperator} using CIL and SIRF took about 350 milliseconds. A comparable implementation in jax \cite{jax2018github} on a T4 GPU from Nvidia took only 2.5 milliseconds.

\section{Conclusions}
Dynamic MRI (like cardiac MRI) is an important clinical imaging modality for which measurement undersampling is vital to get reasonably fast acquisitions. To remedy the image reconstruction artifacts naturally appearing due to this undersampling, motion models such as the optical flow equation can be used to propagate information across time frames. Current models have followed this approach for real-valued images. In this work, we have generalised the procedure to complex-valued images by proposing a complex version of the optical flow equation that can account for phase changes between time frames. We tested the new model on two real cardiac datasets and saw reconstruction improvements of the complex optical flow model over the comparison methods.

We used central differences for spatial derivatives in the optical flow equation. In future work, one could analyse which discrete differentiation approach works best in the context of the complex optical flow equation. Related is the question of whether a CFL condition changes when transitioning from the real to the complex optical flow equation.

We also saw that including ``ground truth'' velocities into the model further improved reconstruction quality implying that the currently chosen regularisation for the velocities is not optimal. In future work, it will be interesting to further investigate how to best choose regularisers for the velocities to achieve better performance (such an optimal choice is not even clear for the standard real-valued optical flow motion model).

Overall, we observed that the complex optical flow model is well-suited for cardiac MRI data. Of course, other options are conceivable. For instance, the continuity equation is a typical way of describing motion (especially in three dimensions). Also, a learning-based approach could be vital to adapt the model to the underlying problem so that problem-specific features can be taken into account. It will be interesting to compare different approaches with each other.

\section*{Acknowledgments}
MM would like to express gratitude to Carola-Bibiane Schönlieb for generously hosting a research stay and supporting this work. MM acknowledges support from
the Deutsche Forschungsgemeinschaft (DFG, German Research
Foundation) under Germany’s Excellence Strategy – EXC 2044 –, Mathematics Münster: Dynamics – Geometry – Structure, and under the Collaborative Research Centre 1450–431460824, InSight, University of Münster, as well as support from the Deutscher Akademischer Austauschdienst (DAAD).

MJE acknowledges support from the EPSRC (EP/S026045/1, EP/T026693/1, EP/V026259/1, EP/Y037286/1) and the European Union Horizon 2020 research and innovation programme under the Marie Skodowska-Curie grant agreement REMODEL.

Calculations for this publication were performed on the HPC cluster PALMA II of the University of Münster, subsidised by the DFG (INST 211/667-1).

Special thanks to Christoph Kolbitsch and Kris Thielemans for their insightful discussions on dynamic MR, which greatly enriched the research process.

%%%%%%%%%%%%%%%%%%%%%%%%%%%%%%%%%%%%%%%%%%%%%%%%%%%%%%
%          7. REFERENCES SECTION
%%%%%%%%%%%%%%%%%%%%%%%%%%%%%%%%%%%%%%%%%%%%%%%%%%%%%%

%       READ THIS SECTION CAREFULLY

% Each of the references below MUST be cited in your article above. Do not include references that are not cited in your article.

% Follow the examples below carefully. We strongly suggest that you copy and paste your reference information directly into our examples.

% List all references in alphabetical order according to the first author's last name.

% Verify each URL works correctly and can be accessed properly. Your URL links should be to reputable websites. The command line for a website link begins with: \url{ }

% Do not add MR or DOI numbers to your references. AIMS production staff will add this information.

% Using BibTex is not recommended but can be handled.

\bibliographystyle{plain}
\bibliography{bibliography}

\end{document}